\documentclass[11pt]{elsarticle}
\usepackage{amsmath,amssymb,amsthm}
\usepackage{hyperref}
\usepackage{mathrsfs}
\usepackage{graphicx}
\usepackage{tikz}

\theoremstyle{plain}
\newtheorem{theorem}{Theorem}[section]
\newtheorem{lemma}[theorem]{Lemma}

\newtheorem{proposition}[theorem]{Proposition}
\newtheorem{corollary}[theorem]{Corollary}

\theoremstyle{definition}
\newtheorem{definition}[theorem]{Definition}
\newtheorem{remark}[theorem]{Remark}
\newtheorem{example}[theorem]{Example}
\newtheorem{problem}[theorem]{Problem}

\newcommand{\ua}{\mathord{\uparrow}}
\newcommand{\da}{\mathord{\downarrow}}
\newcommand{\rom}[1]{\rm{\uppercase\expandafter{\romannumeral #1}}}

\makeatletter
\def\ps@pprintTitle{%
  \let\@oddhead\@empty
  \let\@evenhead\@empty
  \def\@oddfoot{\reset@font\hfil\thepage\hfil}
  \let\@evenfoot\@oddfoot
}
\makeatother

\begin{document}

\begin{frontmatter}

\title{$D$-completion, well-filterification and  sobrification \tnoteref{t1}}
\tnotetext[t1]{This work is supported by the National Natural Science Foundation of China (No.11771134) and by Hunan Provincial Innovation Foundation For Postgraduate (CX20200419)}
\author{Hualin Miao}
\ead{miaohualinmiao@163.com}
\author{Longchun Wang\corref{}}
\ead{longchunw@163.com}
\author{Qingguo Li\corref{a1}}
\address{School of Mathematics, Hunan University, Changsha, Hunan, 410082, China}
\cortext[a1]{Corresponding author.}
\ead{liqingguoli@aliyun.com}
\begin{abstract}
 In this paper, we obtain some sufficient conditions for the $D$-completion of a $T_{0}$ space to be the well-filterification of this space, the well-filterification of a $T_{0}$ space to be the sobrification of this space and the $D$-completion of a $T_{0}$ space to be the sobrification, respectively. Moreover, we give an example to show that a tapered closed set may be neither the closure of a directed set nor the closed $KF$-set, respectively. Because the tapered closed set is a closed $WD$-set, the example also gives a negative answer to a problem proposed by Xu. Meantime, a new direct characterization of the $D$-completion of a $T_{0}$ space is given by using the notion of pre-$c$-compact elements.
\end{abstract}

\end{frontmatter}

\section{Introduction}

 $D$-completion, well-filterification and sobrification of a $T_{0}$ space play a fundamental role in non-Hausdorff topological spaces. We know that the $D$-completion of a $T_{0}$ space is contained in the well-filterification of this space, and the well-filterification is contained in the sobrification. But the converses of them are not necessarily true.

 In \cite{Kei03}, Keimel and Lawson verified that the $D$-completion of a $T_{0}$ space is the sobrification of this space if the space is a $c$-space or $qc$-space. Lawson, Wu and Xi gave the sufficient condition that the space $X$ is core-compact for a well-filtered space $X$ to be sober in \cite{J}. Xu and Shen proved that every first countable well-filtered space is sober in \cite{FX}. Xi and Lawson obtained the result that every monotone convergence space $X$ with the property that $\da (K\cap A)$ is a closed subset of $X$ for any closed subset $A$ and compact saturated set $K$ is well-filtered in \cite{Xi03}. Therefore, naturally, there are some questions in the following:

 (1) Whether the well-filterification of a first countable $T_{0}$ space $X$ coincides with the sobrification of this space;

 (2) Whether the well-filterification of a second countable $T_{0}$ space agrees with the sobrification;

 (3) Whether the $D$-completion of a $T_{0}$ space $X$ with the property that $\da (K\cap A)$ is a closed subset of $X$ for any closed $KF$ set $A$ and compact saturated set $K$ coincides with the well-filterification.

   In \cite{M}, Miao and Li proposed the concept of join-continuous poset and showed that every core-compact join continuous dcpo is sober. They also obtain the result that a monotone convergence space $X$ with the property that $\da(A\cap W)$ is a closed subset of $X$ for any irreducible closed subset $A$ of $X$ and upper set $W$ is a sober space. Naturally, some questions are raised below.

 (4) Whether the $D$-completion of a locally compact $T_{0}$ space $X$ with the property that $\da (K\cap A) $ is a closed subset of $X$ for any closed $KF$ set $A$ and compact saturated set $K$ agrees with the sobrification;

 (5) Whether the $D$-completion of a core-compact join-continuous poset $P$ coincides with the sobrification;

 (6) Whether the $D$-completion of a $T_{0}$ space $X$ with the property that $\da(A\cap W)$ is a closed subset of $X$ for any irreducible closed subset $A$ of $X$ and upper set $W$ coincides with the sobrification.

In \cite{Liu}, Liu, Li and Wu give a counterexample to show that $\mathbf{WD}(X)$ may not agree with $\mathbf{KF}(X)$ for any $T_{0}$ space $X$, which solved the open problem proposed by Xu in \cite{Xu}. Zhang and Li provided a direct characterization of the $D$-completion of a $T_{0}$ space using the tapered closed subsets (\cite{Kou03}). We know that every directed set is a tapered set. So, naturally there is a problem in the following.

\begin{problem}\label{p}
 Let $X$ be a $T_{0}$ space and $A$ a tapered closed subset of $X$. Does there exist a directed subset $D$ of $X$ such that $A=cl(D)$?
\end{problem}
In~\cite{Dcpo03}, Zhao and Fan gave the $D$-completion of a poset. Keimel and Lawson presented a standard $D$-completion of a $T_0$ topological space in~\cite{Kei03}. However, both of their work do not give the concrete elements of  the $D$-completion of a poset or a $T_0$ space like the sobrification of the spaces. In~\cite{Kou03}, Zhang and Li provided a direct description of the elements of the $D$-completion of a $T_{0}$ space by continuous functions. Shall we just characterize the $D$-completion of a $T_0$ space by the elements of the space itself?

 The purpose of this paper is to investigate the above questions. We give some sufficient conditions for the $D$-completion of a $T_{0}$ space to be the well-filterification, the well-filterification of a $T_{0}$ space to be the sobrification and the $D$-completion of a $T_{0}$ space to be the sobrification, respectively. For Problem \ref{p}, we propose a counterexample to reveal that a tapered closed set may not be the closure of a directed set and the closed $KF$-set, respectively. Meantime, a new direct characterization of the $D$-completion of a $T_{0}$ space is given by using the notion of pre-$c$-compact elements.

\section{Preliminaries}
Without further references, the posets mentioned here are all endowed with the Scott topology.

  Let $X$ be a topological space. We denote the set of all \emph{closed sets} of $X$ by $\Gamma (X)$ and all \emph{open sets} by $\mathcal{O} (X)$. If $P$ is a poset, $\sigma (P)$ denotes the set of all \emph{Scott open} sets and $\Gamma(P)$ the set of all \emph{Scott closed} sets. $A$ is called \emph{$d$-closed} if $D\subseteq A$ implies that $\sup D\in A$ for any directed subset of $A$. $cl_{d}(A)$ represents $A$ is $d$-closed. Given a topological space $(X, \tau)$, we define $x\leq y$ iff $x\in cl(y)$. Hence, $X$ with its specialization order is a poset. For any set $ A\subseteq X$, we denote the closure of $A$ by $cl(A)$. The irreducible sets of $X$ are denoted by $IRR(X)$ and the irreducible closed sets by $\mathbf{IRR}(X)$, all upper sets of $X$ by $up(X)$.

\begin{definition}(\cite{clad3})

(\romannumeral1) Let $P$ be a poset. A subset $D$ of $P$ is \emph{directed} provided it is nonempty and every finite subset of $D$ has an upper bound in $D$. We denote that $D$ is a directed subset of $P$ by $D\subseteq ^{\uparrow}P$.

(\romannumeral2) A poset $P$ is a \emph{dcpo} if every directed subset $D$ has the supremum.

\end{definition}

\begin{definition}(\cite{nht2})

(\romannumeral1) A topological space $X$ is \emph{locally compact} if for every $x\in X$ and any open neighborhood $U$ of $x$, there is a compact saturated subset $Q$ of $X$ such that $x\in \mathrm{int}(Q)$ and $Q\subseteq U$.

A set $K$ of a topological space is called \emph{saturated} if it is the intersection of its open neighborhood ($K=\ua K$ in its specialization order). If $A$ is any subset of $X$, the intersection sat$A$ of all its open neighborhood is a saturated set called its \emph{saturation}.

(\romannumeral2) A topological space $X$ is \emph{core-compact} if $\mathcal{O}(X)$ is a continuous lattice.

(\romannumeral3) A topological space $X$ is said to be \emph{first countable} if each point has a countable neighbourhood basis (local base). That is, for each point $x$ in $X$ there exists a sequence $N_{1},N_{2},\cdots$ of the neighbourhoods of $x$ such that for any neighbourhood $N$ of $x$, there exists an integer $i$ with $N_{i}$ contained in $N$.

(\romannumeral4) A topological space is said to be \emph{second countable} if its topology has a countable base.

\end{definition}

\begin{definition}(\cite{nht2})
  A topological space $X$ is \emph{sober} if it is $T_{0}$ and every irreducible closed subset of $X$ is the closure of a (unique) point.
\end{definition}
\begin{definition}(\cite{clad3})
  We shall say that a space $X$ is \emph{well-filtered} if for each filter basis $\mathcal{C}$ of compact saturated sets and each open set $U$ with $\bigcap \mathcal{C} \subseteq U$, there is a $K\in \mathcal{C}$ with $K\subseteq U$.
\end{definition}

 For a topological space $X$, the compact saturated subsets of $X$ are denoted by $Q(X)$. We write $\mathfrak{K} \subseteq_{flt}Q(X)$ represents that $\mathfrak{K}$ is a filtered subfamily of $Q(X)$ and $F\subseteq _{fin}X$ represents that $F$ is a finite subset of $X$.

 \begin{definition}(\cite{She04})\label{a}
   Let $X$ be a topological space. A nonempty subset $A\in KF(X)$ if and only if there exists $ \mathfrak{K} \subseteq_{flt}Q(X)$ such that $cl(A)$ is a minimal closed set that intersects all members of $\mathfrak{K}$. The set of all closed $KF$-subsets of $X$ is denoted by $\mathbf{KF}(X)$.
  \end{definition}

\begin{definition}\cite{FX}
  A $T_{0}$ space $X$ is called \emph{$\omega$-well-filtered}, if for any countable filtered family $\{K_{i} : i < \omega \}\subseteq Q(X)$ and $U \in \mathcal{O}(X)$, it satisfies
  \begin{center}
    $\bigcap_{i<\omega}K_{i}\subseteq U\Rightarrow \exists i_{0}< \omega, K_{i_{0}}\subseteq U$
  \end{center}
\end{definition}
\begin{definition}\cite{FX}
Let $X$ be a $T_{0}$ space. A nonempty subset $A\in KF_{\omega}(X)$ if and only if there exists a countable filtered family $\mathcal {K}\subseteq Q(X)$ such that $cl(A)$ is a minimal closed set that intersects all members of $\mathcal{K}$. The set of all closed $KF_{\omega}$-subsets of $X$ is denoted by $\mathbf{KF_{\omega}}(X)$.
\end{definition}
\begin{definition}\cite{FX}
A subset $A$ of a $T_{0}$ space $X$ is called a \emph{$\omega$-well-filtered determined set}, $WD_{\omega}$ set for short, if for any continuous mapping $f : X\longrightarrow Y$ to a $\omega$-well-filtered space $Y$ , there exists a unique $y_{A} \in Y$ such that $cl(f(A)) = cl(\{y_{A}\})$. The set of all closed $\omega$-well-filtered determined subsets of $X$ is denoted by $\mathbf{WD_{\omega}}(X)$.
\end{definition}

 \begin{definition}(\cite{Kei03})
  A $T_0$ space is a \emph{monotone convergence space} if and only if the closure of every directed sets (in the specialization order) is the closure of a unique point.
\end{definition}

  \begin{definition}(\cite{Kei03})
    A \emph{D-completion} of a $T_0$ space $X$ is a monotone convergence space $Y$ with a topological embedding $j: X \longrightarrow Y $ and $cl_{d}(j(X))=Y$.
   \end{definition}
   \begin{definition}\cite{M}
  A poset $L$ is \emph{join-continuous} if $f_{x}$ is continuous for any $x\in L$ where $f_{x}: (L,\sigma(L))\rightarrow up(L)$ is defined by $f_{x}(y)=\ua x\cap \ua y $ for all $y\in L$ when $up(L)$ is endowed with the upper Vietoris topology.
\end{definition}
\begin{definition}\cite{aou03}
A subset $A$ of a $T_{0}$ space $X$ is called a \emph{well-filtered determined set}, $WD$ set for short, if for any continuous mapping $f : X\longrightarrow Y$ to a well-filtered space $Y$ , there exists a unique $y_{A} \in Y$ such that $cl(f(A)) = cl(\{y_{A}\})$. The set of all closed well-filtered determined subsets of $X$ is denoted by $\mathbf{WD}(X)$.
\end{definition}
\begin{lemma}(\cite{nht2})\label{j}
  If a topological space $X$ is second countable, then $X$ is first countable.
\end{lemma}

 \begin{lemma}\label{jjj}(\cite{Kou03})
The topological space of all tapered closed subsets of a $T_{0}$ space $X$ is the standard $D$-completion of $X$.
 \end{lemma}
 The following construction is due to Ershov \cite{YE}. Let topological spaces $X$ and $Y_{x}$, and $x \in X$, be given.
Let

\begin{center}
  $Z=\bigcup_{x\in X}Y_{x}\times \{x\}$
\end{center}
\begin{center}
  $\tau=\{U\subseteq Z\mid (U)_{x}\in \tau({Y_{x}}) ~for ~any ~x\in X~and~(U)_{X}\in \tau(X) \}$,
\end{center}
where $U_{x}=\{y\in Y_{x}\mid(y,x)\in U\}$ for any $x\in X$ and $(U)_{X}=\{x\in  X\mid (U)_{x}\neq \emptyset\}.
$

In this paper, the space $Z=(Z,\tau)$ is denoted by $\sum_{X}Y_{x}$. For any subset $A \subseteq Z$, put $(A)_{x} = \{y \in Y _{x} \mid
(y,x) \in A\}$, and $(A)_{X} = \{x \in X \mid (A)_{x}\neq \emptyset \}$.
\begin{lemma}\cite{YE}\label{n}
 Let $X$ be a $T_{0}$ space, $Y_{x}$ an irreducible $T_{0}$ space for any $x\in X$, and $Z=\sum_{X}Y_{x}$. For all $(y_{0},x_{0}),(y_{1},x_{1})\in Z$, we have $(y_{0},x_{0})\leq_{Z}(y_{1},x_{1})$ if and only if the following two alternatives hold:

 $(1)$ $x_{0}=x_{1}$, $y_{0}\leq_{Y_{x_{0}}}y_{1}$;

 $(2)$ $x_{0}< x_{1}$ and $y_{1}=\top_{x_{1}}$ is the greatest element in $Y_{x_{1}}$ with respect to the specialization order.
\end{lemma}
\begin{lemma}\cite{Liu}\label{m}
  Let $\mathbb{N} = (N, \tau_{ cof} )$, where $\tau_{cof}$ denotes the cofinite topology, and $X _{n}$ an irreducible $T_{0}$ space for every $n \in N$, such that there are at most
finitely many $X _{n}$'s that have a greatest element under the specialization order $\leq_{X _{n}}$. Let $Z=\sum_{\mathbb{N}} X_{n}$. Then $A$ is a $KF$-set of $Z$ iff there exists a unique $n \in \mathbb{N}$ such that $A \subseteq X_{ n}\times \{n\}$ and $\{y \in X_ {n} \mid (y,x) \in A\}$
is a $KF$-set of $X _{n}$.
\end{lemma}
 \begin{proposition}(\cite{Kei03})\label{iii}
   For a monotone convergence space $X$. $A$ is a monotone convergence subspace iff $A$ is a sub-dcpo of $X$ with the specialization order.
 \end{proposition}

\begin{theorem}\cite{Xu}\label{b}
  Every locally compact $T_{0}$ space is a Rudin space, that is, $\mathbf{IRR}(X)=\mathbf{KF}(X)$.
\end{theorem}
\begin{theorem}\label{f} \cite{meet5}
 Let $L$ be a dcpo. Then the following statements are
equivalent:

$(1)$ $\sigma(L)$ is a continuous lattice;

$(2)$ For every dcpo or complete lattice $S$, one has $\sigma(S \times L) = \sigma (S) \times \sigma (L)$.
\end{theorem}

\begin{theorem}\cite{aou03}\label{o}
Let $X$ be a $T_{0}$ space. Then $\mathbf{WD}(X)$ with the lower Vietoris topology is the well-filtered reflection of $X$.
\end{theorem}

\section{Sufficient conditions}
\begin{lemma}\label{c}
  Let $X$ be a $T_{0}$ space. Then the following statements are equivalent:

  $(1)$ If $A\in \mathbf{KF}(X)$, then $\da (A\cap K)\in \Gamma(X)$ for any $K\in Q(X)$;

  $(2)$ $A$ is a directed closed subset of $X$.
  \begin{proof}
    $(1)\Rightarrow(2)$ It suffices to prove that $A$ is directed. Suppose $A\in \mathbf{KF}(X)$. Then there exists $\{K_{i}\}_{i\in I}\subseteq _{flt} Q(X)$ such that $A$ is a minimal closed set that intersects all members of $\{K_{i}\}_{i\in I}$ by Definition \ref{a}.

    Claim 1: $A=\da (A\cap K_{j})$ for any $j\in I$.

    For any $i\in I$, since $\{K_{i}\}_{i\in I}\subseteq _{flt} Q(X)$, there exists $r\in I$ such that $K_{r}\subseteq K_{i}\cap K_{j}$. This implies that $\emptyset \neq A\cap K_{r}\subseteq A\cap K_{j} \cap K_{i}\subseteq \da(A\cap K_{j})\cap K_{i}$. Then we have $A=\da (A\cap K_{j})$ by the minimality of $A$.

    Claim 2: $\ua x \cap K_{i}\cap A \neq \emptyset$ for any $x\in A$ and $i\in I$.

    $x\in A=\da (A\cap K_{i})$ by Claim 1. Thus $\ua x \cap A \cap K_{i} \neq \emptyset$.

    Claim 3: $A$ is directed.

   Let $x,y\in A$. Then we have $\da(\ua x\cap A)\cap K_{i}\neq \emptyset$ by Claim 2. It follows that $\da (\ua x\cap A)=A$ by the minimality of $A$. Note that $y\in A=\da(\ua x\cap A)$, thus $\ua y\cap \ua x\cap A\neq \emptyset$.

   $(2)\Rightarrow(1)$ It remains to prove that $\da(A\cap K)\in \Gamma(X)$.

   Claim 1: $\da (\ua x \cap A)\in \Gamma(X)$ for any $x\in X$.

   If $\ua x \cap A = \emptyset$, then we have $\da(\ua x \cap A)=\emptyset \in \Gamma(X)$. Else $\ua x \cap A\neq \emptyset$. Then $x\in A$ because $A$ is a lower set. Obviously, $\da (\ua x\cap A)\subseteq A$. Conversely, $\ua a\cap \ua x \cap A \neq \emptyset$ for any $a\in A$ since $A$ is directed. So we have $\da(\ua x\cap A)=A\in \Gamma(X)$.

   Claim 2: $\da (K\cap A)\in \Gamma(X)$.

   If $K\cap A=\emptyset$, then $\da (K\cap A)=\emptyset \in \Gamma(X)$.
   Else $K\cap A\neq \emptyset$. Pick $k\in K\cap A$. This means that $A\subseteq \da(\ua k\cap A)\subseteq \da (K\cap A)\subseteq A$ by Claim 1. Hence, $\da (K\cap A)=A \in \Gamma(X)$.
  \end{proof}
\end{lemma}
\begin{theorem}\label{h}
  Let $X$ be a locally compact $T_{0}$ space with the property $\da(K\cap A)\in \Gamma(X)$ for any $A\in \mathbf{KF}(X)$ and $K\in Q(X)$. Then the $D$-completion of $X$ coincides with the soberification of $X$.
  \begin{proof}
    Let $A\in \mathbf{IRR}(X)$. It suffices to prove that $A$ is directed. By Theorem \ref{b}, we have $\mathbf{IRR}(X)=\mathbf{KF}(X)$. This implies that $A$ is directed by Lemma \ref{c}.
  \end{proof}
\end{theorem}

 Theorem \ref{h} gives a positive answer for the problem $(4)$ in the introduction.
\begin{theorem}\label{e}
  Let $X$ be a $T_{0}$ space with the property $\da(K\cap A)\in \Gamma(X)$ for any $A\in \mathbf{KF}(X)$ and $K\in Q(X)$. If $\mathbf{KF}(X)$ endowed with the lower Vietoris topology is the well-filterification of $X$. Then the $D$-completion of $X$ agrees with the well-filterification of $X$.
  \begin{proof}
    By Lemma \ref{jjj}, we know that all tapered closed subsets of $X$ is the $D$-completion of $X$. Note that directed closed subsets of $X$ must be tapered. We need to prove that $A$ is a directed subset of $X$ for any tapered closed subset $A$ of $X$. It suffices to prove that $\mathcal{D}=\{cl(D)\mid D\subseteq ^{\uparrow} X\}$ is a subdcpo of $\Gamma(X)$. Let $\{cl(D_{i})\}_{i\in I}$ be any directed subset of $\Gamma(X)$ contained in $\mathcal{D}$. Note that $\sup_{i\in I}cl(D_{i})=cl(\bigcup_{i\in I}cl(D_{i}))$. We need to verify that $\sup_{i\in I}cl(D_{i})\in \mathcal{D}$. Since $\mathcal{D}\subseteq \mathbf{KF}(X)$, we have $cl(D_{i})$ is directed for any $i\in I$ by Lemma \ref{c}. It suffices to show that $\bigcup_{i\in I}cl(D_{i})$ is a directed subset of $X$. Now let $x,y \in \bigcup_{i\in I}cl(D_{i})$. Then there exists $\{i_{x},i_{y}\}\subseteq I$ such that $x\in cl(D_{i_{x}}), y\in cl(D_{i_{y}})$. It follows that there exists $i\in I$ such that $\{x,y\}\subseteq cl(D_{i_{x}})\cup cl(D_{i_{y}})\subseteq cl(D_{i})$ by the directionality of $(cl(D_{i}))_{i\in I}$. Thus there exists $z\in cl(D_{i}) \subseteq \bigcup_{i\in I}cl(D_{i})$ such that $z$ is an upper bound of $x,y$ because $cl(D_{i})$ is a directed subset of $X$. So $\sup_{i\in I}cl(D_{i})\in \mathcal{D}$. By Lemma \ref{c}, we have $\mathbf{KF}(X)\subseteq \mathcal{D}$. Hence, the $D$-completion of $X$ agrees with the well-filterification of $X$.
  \end{proof}
\end{theorem}
We know that a monotone convergence space with the property $\da (K\cap A)\in \Gamma(X)$ for any $A\in \mathbf{KF}(X)$ and $K\in Q(X)$ is a well-filtered space by \cite{Xi03}. So there is a problem below:

\begin{problem}\label{aa}
Without the condition that $\mathbf{KF}(X)$ endowed with the lower Vietoris topology is the well-filterification of $X$, whether the statement in Theorem \ref{e} holds for any $T_{0}$ space $X$ with the property $\da (K\cap A)\in \Gamma(X)$ for any $A\in \mathbf{KF}(X)$ and $K\in Q(X)$.
\end{problem}
Note that Problem \ref{aa} actually is the problem $(3)$ in the introduction. For this problem, we give a counterexample as follows.

\begin{example}
  Let $X=(\mathbb{N},\tau_{cof}), Y_{n}=(\mathbb{N},\sigma(\mathbb{N}))$ for any $n\in X$, $Z=\sum_{X}Y_{n}$, where $\mathbb{N}$ is the set of natural numbers. Then $\da(K\cap A)\in \Gamma(Z)$ for any $A\in \mathbf{KF}(Z)$ and $K\in Q(Z)$, but the $D$-completion $Z^{d}$ of $Z$ does not agree with the well-filterification $Z^{w}$ of $Z$.
  \begin{proof}
    By Lemma \ref{m} and Lemma \ref{n}, we have that $A$ is a directed subset of $Z$ for any $A\in \mathbf{KF}(Z)$. This implies that $\da(K\cap A)\in \Gamma(Z)$ for any $A\in \mathbf{KF}(Z)$ and $K\in Q(Z)$ because of Lemma \ref{c}. Let $\mathcal{N}=\{Y_{n}\times \{n\}\mid n\in X\}$. Note that $Y_{n}\times \{n\} \subseteq ^{\ua}Z $ for any $n\in X$. It follows that $\mathcal{N} \subseteq Z^{w}$.

    Claim 1: $\mathcal{A}$ is compact in $Z^{w}$ for any $\mathcal{A} \subseteq \mathcal{N}$.

    Let $\{U_{i}\}_{i\in I} \subseteq \mathcal{O}(Z)$ with $\mathcal{A}\subseteq \bigcup_{i\in I}\diamondsuit U_{i}$. Pick $Y_{n_{0}}\times \{n_{0}\}\in \mathcal{A}\subseteq \bigcup_{i\in I}\diamondsuit U_{i}$. So there exists $i_{n_{0}}\in I$ such that $Y_{n_{0}}\times \{n_{0}\}\in \diamondsuit U_{i_{n_{0}}}$. $(U_{i_{n_{0}}})_{X}\in \tau_{cof}$ since $U_{i_{n_{0}}}$ is open in $Z$. Let $B=X\backslash (U_{i_{n_{0}}})_{X}, C=\{n\in X\mid Y_{n}\times \{n\}\in \mathcal{A}, Y_{n}\times \{n\} \cap U_{i_{n_{0}}}=\emptyset\}$. Note that $B$ is finite and $C\subseteq B$. Then there exists finitely many members of $\{\diamondsuit U_{i}\}_{i\in I}$ to cover $\mathcal {A}$.

    Let $\mathcal{A}_{n}=\mathcal{N}\backslash \{Y_{i}\times\{i\}\mid i\in \{1,2,\cdots n\} \}$. Thus $\{\ua_{Z^{w}}\mathcal{A}_{n}\}_{n\in X}\subseteq_{flt}Q(Z^{w})$ by Claim 1.

Since $Z^{w}$ is well-filtered, we have $\bigcap_{n\in X}\ua_{Z^{w}}\mathcal{A}_{n}\neq \emptyset$. Choose $B\in \bigcap_{n\in X}\ua_{Z^{w}}\mathcal{A}_{n}$, i.e., $B\in Z^{w}$.

    Claim 2: $(B)_{X}$ is infinite.

    Suppose not, if $(B)_{X}$ is finite. Let $n_{1}=\max{(B)_{X}}$. Since $B\in \bigcap_{n\in X}\ua_{Z^{w}}\mathcal{A}_{n}\subseteq \ua _{Z^{w}} \mathcal{A}_{n_{1}}$, we have $\da_{Z^{w}}B\cap \mathcal{A}_{n_{1}}\neq \emptyset $. Then there exists $n_{2}>n_{1}$ such that $Y_{n_{2}}\times \{n_{2}\}\in \da_{Z^{w}}B\cap \mathcal{A}_{n_{1}}$, which contradicts $n_{2}\notin (B)_{X} $.

    By Lemma \ref{n}, $B$ is not directed. This implies that $B$ is not a tapered closed subset of $Z$ by the proof of Theorem \ref{e}. Thus $B\in Z^{w}\backslash Z^{d}$, that is, the $D$-completion $Z^{d}$ of $Z$ does not agree with the well-filterification $Z^{w}$ of $Z$.

  \end{proof}
\end{example}
\begin{lemma}\label{t}
 Let $L$ be a poset. Then the following statements are
equivalent:

$(1)$ $\sigma(L)$ is a continuous lattice;

$(2)$ For every dcpo or complete lattice $S$, one has $\sigma(S \times L) = \sigma (S) \times \sigma (L)$.
\begin{proof}
  The proof is similar to Theorem \ref{f}.
\end{proof}
\end{lemma}
\begin{theorem}\label{g}
  Let $L$ be a core-compact and join continuous poset. Then the $D$-completion of $L$ coincides with the soberification of $L$.
  \begin{proof}
Let $A$ be a irreducible closed subset of $X$. Then it suffices to prove that $A$ is directed. Define $F:(L\times L, \sigma(L\times L))\longrightarrow up(L)$ by $F(x,y)=\ua x\cap \ua y$ for any $(x,y)\in L\times L$.

Claim 1: $F$ is continuous.

Let $U\in \sigma (L)$. We need to prove that $F^{-1}(\Box U)$ is Scott open in $L\times L$. It suffices to prove that $F^{-1}(\Box U)$ is Scott open when $F^{-1}(\Box U)\neq \emptyset$. Obviously, $F^{-1}(\Box U)=\ua (F^{-1}(\Box U))$.
Now let $\{(x_{i},y_{i})\}_{i\in I}$ be any directed subset of $L\times L$ such that $\sup _{i\in I}(x_{i},y_{i})$ exists in $L\times L$ with $\sup_{i\in I}(x_{i},y_{i})=(\sup_{i\in I}x_{i}, \sup_{i\in I}y_{i}) \in F^{-1}(\Box U)$. This implies that $\ua \sup_{i\in I}x_{i}\cap \ua \sup_{i\in I}y_{i}\in \Box U$. So we have $f_{\sup_{i\in I}x_{i}}(\sup_{i\in I}y_{i})\in \Box U$. It follows that $\sup_{i\in I}y_{i}\in f_{\sup_{i\in I}x_{i}}^{-1}(\Box U)$ is Scott open in $L$ because $L$ is join continuous. This implies that there exists $i_{0}\in I$ such that $y_{i_{0}}\in f_{\sup_{i\in I}x_{i}}^{-1}(\Box U)$. Thus $\ua y_{i_{0}}\cap \ua \sup_{i\in I}x_{i}\in \Box U$. Then we have $\sup_{i\in I}x_{i}\in f_{y_{i_{0}}}^{-1}(\Box U)$. By the join continuity of $L$, there exists $i_{1}\in I$ such that $x_{i_{1}}\in f_{y_{i_{0}}}^{-1}(\Box U)$. Thus there exists $i\in I$ such that $(x_{i},y_{i})$ is an upper bound of $\{(x_{i_{0}},y_{i_{0}}),(x_{i_{1}},y_{i_{1}}) \}$ because $\{(x_{i},y_{i})\}_{i\in I}$ is directed. Therefore, $(x_{i},y_{i})\in F^{-1}(\Box U)$. Hence, $F$ is continuous.

Claim 2: $A$ is a directed subset of $X$.

By Lemma \ref{t} and Claim 1, we have $F:\sigma(L)\times \sigma(L)\longrightarrow up(L)$ is continuous. Let $\{x,y\}\subseteq A $. We need to prove that $\ua x \cap \ua y\cap A \neq \emptyset$. Suppose not, if $\ua x \cap \ua y \cap A=\emptyset$. It follows that $\ua x\cap \ua y\subseteq \Box (X\backslash A)$, that is, $F((x,y))\in \Box (X\backslash A)$. This implies that there exists $\{U,V\}\subseteq \sigma(L)$ such that $(x,y)\in U\times V\subseteq F^{-1}(\Box (X\backslash A))$ by the continuity of $F$. Notice that $x\in A\cap U$, $y\in A\cap V$. We have $A\cap U\cap V\neq \emptyset$ since $A$ is irreducible. Pick $z\in A\cap U\cap V$. Hence $(z,z)\in U\times V\subseteq F^{-1}(\Box(X\backslash A))$. This means that $z\in \ua z\cap \ua z\subseteq (X\backslash A)$, which contradicts $z\in A$. So $A$ is directed.

  \end{proof}
\end{theorem}
For the problem $(5)$ mentioned in the introduction, we give a positive answer by Theorem \ref{g}.
\begin{corollary}\label{i}
  Let $L$ be a core-compact poset with the property that $\da(K\cap A)\in \Gamma(L)$ for any $A\in \mathbf{IRR}(L)$ and $K\in Q(L)$. Then the $D$-completion of $L$ agrees with the soberification of $L$.
  \begin{proof}
    The proof is similar to Theorem \ref{g}.
  \end{proof}
\end{corollary}
In contrast to Theorems \ref{h} and Corollary \ref{i}, we raise the following question.
\begin{problem}
    Whether the statement in Corollary \ref{i} holds for arbitrary core-compact topological spaces $X$ with the property that $\da (K\cap A)\in \Gamma(X)$ for any $A\in \mathbf{IRR}(X)$ and $K\in Q(X)$.
\end{problem}
\begin{proposition}\label{w}
  Let $X$ be a second countable space. Then the well-filterification of $X$ coincides with the soberification of $X$.
  \begin{proof}
    It suffices to prove that the well-filterification of $X$ is second countable by Theorem 4.1 in \cite{FX}. Let $\mathcal{B}$ be a countable basis of $X$, $\mathcal{U}=\{\diamondsuit U\mid U\in \mathcal{B}\}$. Since $\mathcal{U}$ is countable, we only need to prove that $\mathcal{U}$ is a basis of the well-filterification of $X$. Let $A\in \mathbf{WD}(X)$, $V\in \mathcal{O}(X)$, $A\in \diamondsuit V$. Then $A\cap V\neq \emptyset$. Pick $x\in A\cap V$. This implies that there exists $U\in \mathcal{B}$ such that $x\in U\subseteq V$. It follows that $A\in \diamondsuit U\subseteq \diamondsuit V$.
  \end{proof}
  We give a positive answer for the problem $(2)$ in the introduction by Theorem \ref{w}.
\end{proposition}
\begin{lemma}\label{v}
  Let $X$ be a $T_{0}$ space, $X^{s}$ the soberification of $X$. For each ordinal $\beta$, define

  $(1)$ $X_{\beta+1}=\{x\in X^{s}\mid \exists F\in \mathbf{KF}(X_{\beta}),cl_{X^{s}}(F)=\da _{X^{s}}x\}$;

  $(2)$ $X_{\beta}=\bigcup_{\gamma< \beta}X_{\gamma}$ for a limit ordinal $\beta$.

  Then there exists an ordinal $\alpha_{X}$ such that $X_{\alpha_{X}}=X_{\alpha_{X}+1}$ and the $\omega$-well-filtered reflection is $X_{\alpha_{X}}$ endowed with the lower Vietoris topology.
  \begin{proof}
    The proof is similar to Proposition 3.8 in \cite{She04}.
  \end{proof}
\end{lemma}
\begin{theorem}\label{l}
  Let $X$ be a first countable space. Then the well-filterification of $X$ agrees with the soberification of $X$.
  \begin{proof}
    By Theorem 4.1 in \cite{FX}, we only need to prove that the $\omega$-well-filterification of $X$ is first countable. There exists an ordinal $\alpha_{X}$ such that $X_{\alpha_{X}}=X_{\alpha_{X}+1}$ and the $\omega$-well-filtered reflection is $X_{\alpha_{X}}$ endowed with the lower Vietoris topology by Lemma \ref{v}. We want to prove that $X_{\alpha}$ endowed with the lower Vietoris topology is first countable for any ordinal $\alpha\leq \alpha_{X}$. We use induction on $\alpha$. For $\alpha=0$, $X_{0}=\{\da x\mid x\in X\}$. Then $X$ and $X_{0}$ are homeomorphic. This implies that $X_{0}$ is first countable. Let $\alpha$ be such that $\alpha+1\leq \alpha_{X}$, and let $X_{\alpha}$ be first countable.

    Claim 1: $X_{\alpha +1}$ is first countable.

    Let $A\in X_{\alpha +1}$. Then there exists $\mathcal{F}\in \mathbf{KF}(X_{\alpha})$ such that $cl_{X^{s}}(\mathcal{F})=\da _{X^{s}}A$. Since $\mathcal{F}\in \mathbf{KF}(X_{\alpha})$, there exists  $(\mathcal {K}_{n})_{n\in \mathbb{N}} \subseteq Q(X_{\alpha})$ such that $\mathcal{F}$ is a minimal closed set that intersects all members of $\{\mathcal {K}_{n}\mid n\in \mathbb{N}\}$. Pick $A_{n}\in \mathcal{F}\cap \mathcal{K}_{n}$ for any $n\in \mathbb{N}$. It follows that $\mathcal{F}=cl_{X_{\alpha}}(\mathcal{A})$ from the minimality of $\mathcal{F}$, where $\mathcal{A}=\{A_{n}\mid n\in \mathcal{N}\}$. We want to prove that $cl_{X^{s}}(\mathcal {F})=cl_{X^{s}}(\mathcal{A})$. Since $cl_{X^{s}}(\mathcal{F})\cap X_{\alpha}=\mathcal{F}= cl_{X_{\alpha}}(\mathcal{A})=cl_{X^{s}}(\mathcal{A})\cap X_{\alpha}$, we have $\mathcal{F}\subseteq cl_{X^{s}}(\mathcal{A})$. So $cl_{X^{s}}(\mathcal{F})\subseteq cl_{X^{s}}(\mathcal{A})$. Conversely, note that $A_{n}\in \mathcal{K}_{n}\subseteq X_{\alpha}$ for any $n\in \mathbb{N}$. Hence, $\mathcal{A}\subseteq cl_{X^{s}}(\mathcal{A})\cap X_{\alpha}=cl_{X^{s}}(\mathcal{F})\cap X_{\alpha}$. This implies that $\mathcal{A}\subseteq cl_{X^{s}}(\mathcal{F})$. Therefore, $cl_{X^{s}}(\mathcal {F})=cl_{X^{s}}(\mathcal{A})$. Because $X_{\alpha}$ is first countable, there exists a countable neighborhood basis $\mathcal{B}_{n}=\{\diamondsuit U\cap X_{\alpha}\mid \diamondsuit U\cap X_{\alpha}\in \mathcal{B}_{n}\}$ of $A_{n}$ for any $n\in \mathbb{N}$. Let $\mathcal{B}_{A}=\{\diamondsuit U\cap X_{\alpha+1}\mid \diamondsuit U\cap X_{\alpha}\in \bigcup_{n\in\mathbb{N}}\mathcal{B}_{n}\}$. Note that $\mathcal{B}_{A}$ is countable. We only need to prove that $\mathcal{B}_{A}$ is a countable neighborhood basis of $A$ in $X_{\alpha+1}$. Let $V\in \mathcal{O}(X)$ with $A\in \diamondsuit V\cap X_{\alpha+1}$. It follows that $\mathcal{A}\cap \diamondsuit V\neq \emptyset$ since $\da_{X^{s}}A=cl_{X^{s}}(\mathcal{F})=cl_{X^{s}}(\mathcal{A})$. Choose $A_{n}\in \mathcal{A}\cap \diamondsuit V\cap X_{\alpha}$. Then there exists $\diamondsuit U\cap X_{\alpha}\in \mathcal{B}_{n}$ such that $A_{n}\in \diamondsuit U\cap X_{\alpha}\subseteq \diamondsuit V\cap X_{\alpha}$. We need to prove that $U\subseteq V$. Suppose not, if $U\nsubseteq V$, then there exists $u\in U\cap (X\backslash V)$. This implies that $\da u\in \diamondsuit U\cap X_{\alpha}\subseteq \diamondsuit V$. Thus $u\in V$, which contradicts $u\in X\backslash V$. So $A \in \diamondsuit U\cap X_{\alpha+1} \subseteq \diamondsuit V\cap X_{\alpha +1}$.

    Suppose now that $\alpha \leq \alpha_{X}$ is a limit ordinal and $X_{\beta}$ is first countable for any $\beta <\alpha$. By definition of $X_{\alpha}$, we have $X_{\alpha}=\bigcup_{\beta<\alpha}X_{\alpha}$.

    Claim 2: $X_{\alpha}$ is first countable.

    For any $A\in X_{\alpha}$, there exists $\beta<\alpha$ such that $A\in X_{\beta}$.
    We have that there exists a countable neighborhood basis $\mathcal{B}=\{\diamondsuit U\cap X_{\beta}\mid \diamondsuit U\cap X_{\beta}\in \mathcal{B}\}$ of $A$ in $X_{\beta}$ since $X_{\beta}$ is first countable. It suffices to prove that $\mathcal{B}_{A}=\{\diamondsuit U\cap X_{\alpha}\mid \diamondsuit U\cap X_{\beta}\in \mathcal{B}\}$ is a countable neighborhood basis of $A$ in $X_{\alpha}$. It is obvious that $\mathcal{B}_{A}$ is countable. Now let $V\in \mathcal{O}(X)$ with $A\in X_{\alpha}\cap \diamondsuit V$. Then $A\in \diamondsuit V \cap X_{\beta} $. It follows that there exists $\diamondsuit U\cap X_{\beta}\in \mathcal{B}$ such that $A\in \diamondsuit U\cap X_{\beta}\subseteq \diamondsuit V\cap X_{\beta}$. Similar to the proof of Claim 1, we have $U\subseteq V$. So $A\in \diamondsuit U\cap X_{\alpha}\subseteq \diamondsuit V\cap X_{\alpha}$. This proves that $X_{\alpha}$ endowed with the lower Vietoris topology is first countable for any ordinal $\alpha\leq \alpha_{X}$. Hence, the $\omega$-well-filtered reflection $X_{\alpha_{X}}$ of $X$ is first countable.

  \end{proof}
\end{theorem}
 We give a positive answer by Theorem \ref{l} to the problem $(1)$ in the introduction.
\begin{theorem}\label{a1}
  Let $X$ be a $T_{0}$ space with the property $\da (A\cap W)\in \Gamma(X)$ for any $A\in \mathbf{IRR}(X)$ and $W\in up(X)$. Then the $D$-completion of $X$ coincides with the soberification of $X$.
  \begin{proof}
  Let $A\in \mathbf{IRR}(X)$. It suffices to prove that $A$ is directed. For any $x,y \in A$, we need to prove that $\ua x\cap \ua y \cap A\neq \emptyset$. Suppose not, $\ua x\cap \ua y \cap A=\emptyset$. Let $B=\da(\ua x \cap A)$. Then $\da (X\backslash B)\in \Gamma(X)$. It follows that $A=(B\cup (X\backslash B))\cap A=(B\cap A)\cup ((X\backslash B)\cap A)=B\cup ((X\backslash B)\cap A)\subseteq B\cup \da((X\backslash B)\cap A)$. Thus $A\subseteq B$ or $A\subseteq \da((X\backslash B)\cap A)$. But $y \in A\cap (X\backslash B) $ and $\ua x\cap A\cap (X\backslash B)=\emptyset$ implies that $A\nsubseteq B$ and $A\nsubseteq \da((X\backslash B)\cap A)$, which contradicts $A\subseteq B$ or $A\subseteq \da((X\backslash B)\cap A)$.
  \end{proof}
\end{theorem}
For the problem $(6)$ in the introduction, we give a positive answer by Theorem \ref{a1}.
\section{A counterexample}

In the following, we construct an example to give a negative answer to Problem \ref{p}. It also can answer Xu's problem (\cite{Xu}).
\begin{figure}[htbp]
\centering
\includegraphics[height=7cm, width=16cm]{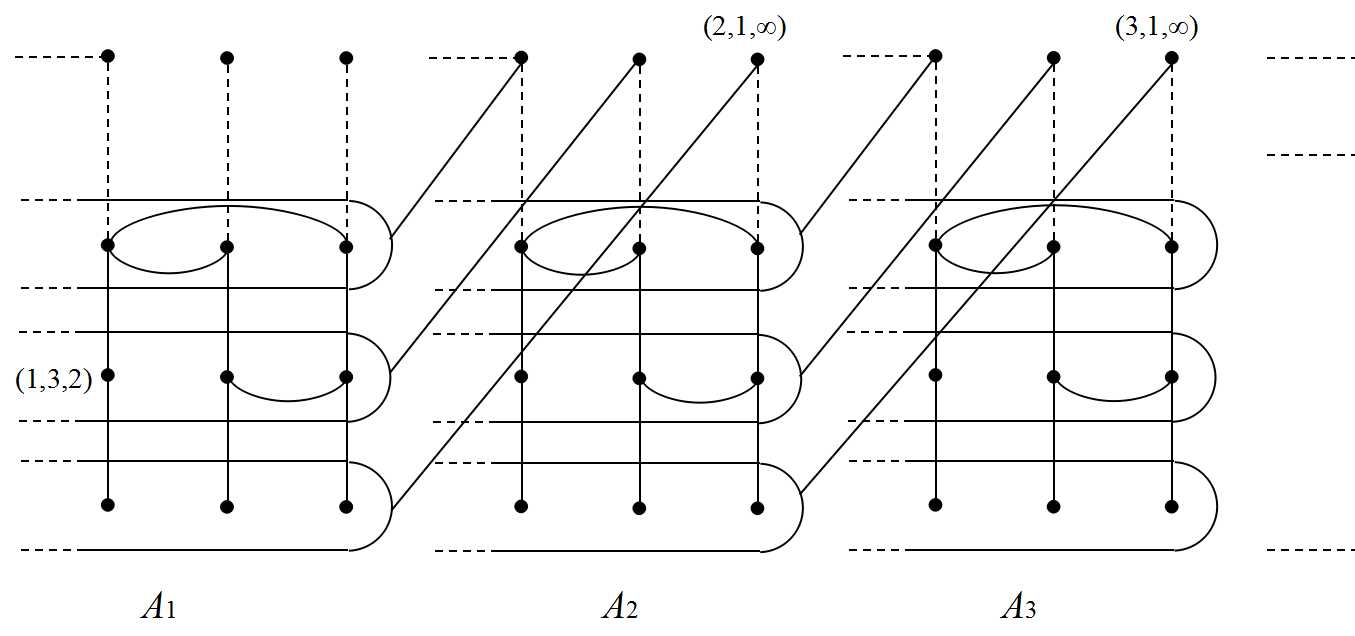}
\caption{A tapered non-directed-determined poset.}\label{r}
\end{figure}
\begin{example}\label{q}
  Let $L=\mathbb {N}\times \mathbb{N}\times (\mathbb{N}\cup \infty)$, where $\mathbb{N}$ is the set of natural numbers. We define an order $\leq$ on $L$ as follows:

  $(n_{1},i_{1},j_{1})\leq (n_{2},i_{2},j_{2})$ if and only if:

  $\bullet$ $n_{1}=n_{2}, i_{1}=i_{2}, j_{1}\leq j_{2}$;

  $\bullet$ $n_{1}=n_{2},i_{2}=j_{2}=j_{1}, i_{1}\leq i_{2}$;

  $\bullet$ $n_{2}=n_{1}+1, j_{1}\leq i_{2},j_{2}=\infty$.

  $L$ can be easily depicted as in Figure \ref{r}. Then $L$ is tapered closed in $(L,\sigma(L))$, but $L\notin \mathbf{KF}(L)$ and $L\notin \{cl(D)\mid D\subseteq^{\uparrow}L\}$.
  \begin{proof}
    Claim 1: $L$ is tapered closed.

   Let $A_{n}=\{(n,i,j)\mid (i,j)\in \mathbb{N}\times (\mathbb{N}\cup \infty)\}$. Note that $\bigcup_{i=1}^{n}A_{i}=cl(A_{n}\backslash \max A_{n})$ for any $n\in \mathbb{N}$ and $A_{n}\backslash \max A_{n}\subseteq^{\uparrow}L$. Hence, $\bigcup_{i=1}^{n}A_{i}$ is tapered closed. Let $B_{n}=\bigcup_{i=1}^{n}A_{i}$. It follows that $\{B_{n}\}_{n\in \mathbb{N}}\subseteq ^{\ua} \Gamma(L)$. Then we have $L=\sup_{n\in \mathbb{N}}B_{n}$ is tapered closed.

   Claim 2: $L\notin \{cl(D)\mid D\subseteq^{\uparrow}L\}$.

   Suppose not, if there exists a directed subset $D$ of $L$ such that $cl(D)=L$. We need to prove that $(D)_{\mathbb{N}}=\{n\in \mathbb{N}\mid A_{n}\cap D \neq \emptyset\}$ is infinite. Assume $(D)_{\mathbb{N}}$ is finite. Let $n_{0}=\max(D)_{\mathbb{N}}$. Then we have $D\subseteq\bigcup_{i\in (D)_{\mathbb{N}}}A_{i}\subseteq \bigcup_{i=1}^{n_{0}}A_{i}=B_{n_{0}}$. So $L= cl(D)\subseteq B_{n_{0}}$, which contradicts $B_{n_{0}}\subsetneqq L$. Thus $(D)_{\mathbb{N}}$ is infinite. Let $n_{1}\in (D)_{\mathbb{N}}$. Pick $(n_{1},i_{1},j_{1})\in D\cap A_{n_{1}}$. By the infiniteness of $(D)_{\mathbb{N}}$, there exists $m\in (D)_{\mathbb{N}}$ such that $m\geq n_{1}+2$. Choose $(m,i_{m},j_{m})\in D \cap A_{m}$. Obviously, $\ua (n_{1},i_{1},j_{1})\subseteq A_{n_{1}}\cup A_{n_{1}+1} $ and $\ua (m,i_{m},j_{m})\subseteq A_{m}\cup A_{m+1}$, $(A_{n_{1}}\cup A_{n_{1}+1})\cap(A_{m}\cup A_{m+1})=\emptyset$. It follows that $\ua (n_{1},i_{1},j_{1})\cap \ua (m,i_{m},j_{m})=\emptyset$, which contradicts that $D$ is directed.

   Claim 3: $(K)_{\mathbb{N}}=\{n\in \mathbb{N}\mid K\cap A_{n}\neq \emptyset\}$ is finite for any $K\in Q(L)$.

   Suppose not, if $(K)_{\mathbb{N}}$ is infinite. For any $n\in (K)_{\mathbb{N}}$, pick $(n,i_{n},j)\in K\cap A_{n}$, Then $(n,i_{n},\infty)\in K\cap A_{n}$. Let $F\subseteq_{fin}(K)_{\mathbb{N}}$, $B_{F}=\{\bigcup_{n\in (K)_{\mathbb{N}}\backslash F} \da(n,i_{n},\infty)\}$. We need to verify that $B_{F}\in \Gamma(L)$. Now let $D\subseteq^{\ua}B_{F}$. If $D$ is finite, then $\sup D\in D\subseteq B_{F}$. Else, there exists a unique $n\in (K)_{\mathbb{N}}$ such that $D\subseteq \da(n,i_{n},\infty)$. This implies that $\sup D\in  \da(n,i_{n},\infty)\subseteq B_{F}$. $\da B_{F}=B_{F}$ is obvious. Note that $B_{F}\cap K\neq \emptyset$ for any $F\subseteq _{fin} (K)_{\mathbb{N}}$. So we have $K\cap \bigcap_{F\subseteq_{fin} (K)_{\mathbb{N}}}B_{F}\neq \emptyset$. Therefore, there exists $(n,i,j)\in K\cap \bigcap_{F\subseteq_{fin} (K)_{\mathbb{N}}}B_{F}$. Since $(K)_{\mathbb{N}}$ is infinite, we have $m\in (K)_{\mathbb{N}}$ such that $m\geq n+2$.

   We now distinguish two cases:

   Case 1, $n+1\in (K)_{\mathbb{N}}$: Then $(n,i,j)\in \bigcap_{F\subseteq_{fin} (K)_{\mathbb{N}}}B_{F}\subseteq B_{\{n,n+1\}}\cap A_{n}$, which contradicts $B_{\{n,n+1\}}\cap A_{n}=\emptyset$.

   Case 2, $n+1\notin (K)_{\mathbb{N}}$: Then $(n,i,j)\in \bigcap_{F\subseteq_{fin} (K)_{\mathbb{N}}}B_{F}\subseteq B_{\{n\}}\cap A_{n}$, which contradicts $B_{\{n\}}\cap A_{n}=\emptyset$.

   Claim 4: $L\notin \mathbf{KF}(L)$.

   Suppose not, $L\in \mathbf{KF}(L)$, then there exists $\{K_{i}\}_{i\in I} \subseteq_{flt}Q(L)$ such that L is a minimal closed set that intersects all members of $\{K_{i}\}_{i\in I}$ by Definition \ref{a}. Pick $i_{0}\in I$. $(K_{i_{0}})_{\mathbb{N}}$ is finite by Claim 3. Now let $m_{0}=\max(K_{i_{0}})_{\mathbb{N}}$. This implies that $K_{i_{0}}\subseteq B_{m_{0}}$. We need to prove that $B_{m_{0}}\cap K_{i}\neq \emptyset$ for any $i\in I$. For any $i\in I$, there exists $j\in I$ such that $K_{j}\subseteq K_{i}\cap K_{i_{0}}$ because $\{K_{i}\}_{i\in I}$ is filter. It follows that $\emptyset \neq K_{j}=B_{m_{0}}\cap K_{j}\subseteq B_{m_{0}} \cap K_{i}$. So $L=B_{m_{0}}$ by the minimality of $L$, which contradicts $L\neq B_{m_{0}}$.

     \end{proof}
\begin{figure}[htbp]
\centering
\includegraphics[height=5cm, width=16cm]{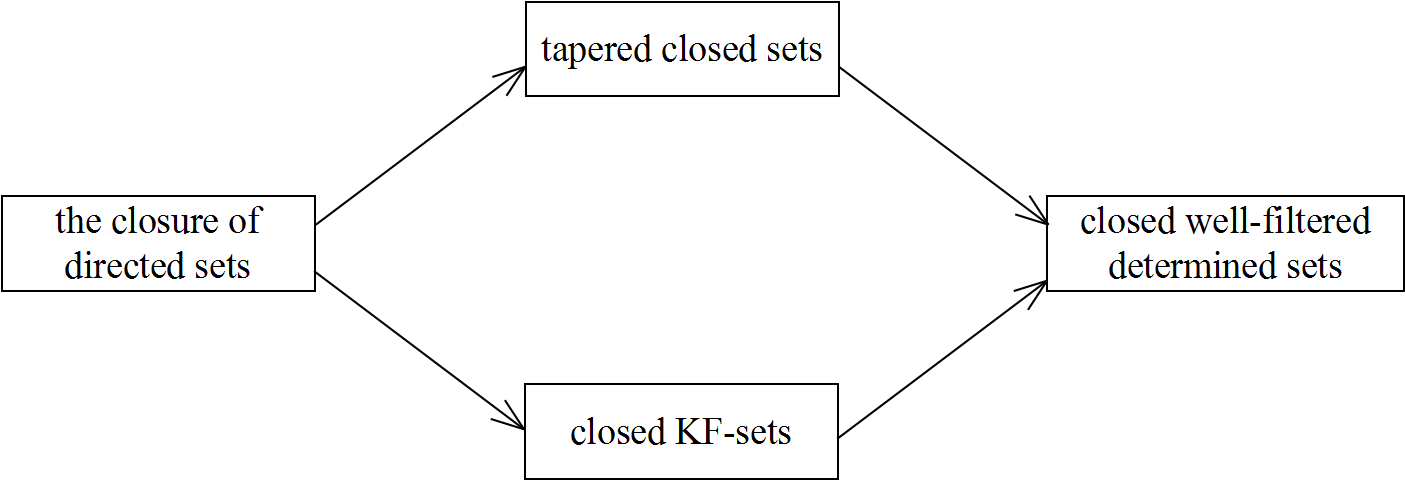}
\caption{Certain relations among some kinds of subsets of a $T_{0}$ space.}\label{u}
\end{figure}
\end{example}
\begin{example}\label{s}
  Let $J$ be Johnstone's Example. Then $J\in \mathbf{KF}(J)$, but $J$ is not a tapered closed set.
\end{example}
\begin{remark}
  By Example \ref{q} and \ref{s}, we know that there is no relation between the tapered closed subsets and the closed $KF$-subsets of a $T_{0}$ space $X$.
\end{remark}
Figure \ref{u} below shows certain relations among some kinds of subsets of a $T_{0}$ space.
\section{A direct characterization for $D$-completion}
    Let $P$ be a poset, $\mathcal{H}=\{\mathcal{D} \subseteq \Gamma(P)\backslash \emptyset \mid  ,\forall A \in \mathcal{D}, \forall a \in A, \da a \in \mathcal{D}, cl_{d}(\mathcal{D})=\mathcal{D}\}$ and $A,B \in \Gamma(P) $. We say that $A ~ \mathrm{is} ~\mathrm{beneath}~B$ denoted by $A\prec B$, if for every $\mathcal{D} \in \mathcal{H} $, the relation $B\subseteq \sup \mathcal{D}$ always implies that $A \in \mathcal{D}$.

Notice that the relation which we define above is not an auxiliary relation because the elements of $\mathcal{H}$ may not be lower sets.
\begin{definition}
  Let $P$ be a poset. An element $A$ of $\Gamma (P)$ is called \emph{pre-C-compact} if $A\prec A$. We use $\mathcal{K} (\Gamma(P))$ to denote the set of all the pre-$C$-compact elements of $P$.
\end{definition}

Ho and Zhao introduce \emph{C-continuous} posets and $C$-compact elements in~\cite{Ho03}. They proved that the $C$-compact elements of a poset $P$ are co-primes. They find  that the $C$-compact elements of $\Gamma(P)$ are more than the elements of the $D$-completion of the posets. Now we want to find a way to reduce the $C$-compact elements of $P$ such that they are equal.

\begin{definition}\cite{Ho03}\label{a2}
   Let $P$ be a poset and $x,y \in P$. We say that $x$ is beneath $y$,
denoted by $x \prec y$, if for every nonempty Scott closed set $C \subseteq P$ for which $\sup C$ exists, the relation $y\leq \sup C$ always implies that $x \in C$. An element $x$ of a poset $P$ is called $C$-compact if $x \prec x$.
\end{definition}

 The pre-$C$-compact elements are $C$-compact by Definition \ref{a2}. The following example reveals that $C$-compact elements may not be pre-$C$-compact.
\begin{example}
Let $J=\mathbb{N} \times (\mathbb{N}\cup \{\infty\})$, that is Johnstone's example. The partial order is $(m,n)\leq (a,b)~\mathrm{iff}~ \mathrm{either}~ m=a ~\mathrm{and} ~n\leq b ~\mathrm{or} ~b=\infty ~\mathrm{and} ~n\leq a$. That $J$ is $C$-compact have been proved in Remark 5.3 of~\cite{Ho03}. We now show that $J$ is not pre-$C$-compact. Let $\mathcal{D}=\{\da x\mid x\in J\}$. Note that $\mathcal{D} \in \mathcal{H}$ and $\sup \mathcal{D}=J$ but $J\not\in \mathcal{D}$.
\end{example}
The following proposition reveals some order-theoretic properties of the lattice $(\Gamma(P),\subseteq)$ for an arbitrary poset $P$.
\begin{proposition}\label{aaa}
   Let P be a poset. Then for any $\mathcal{D} \in \mathcal{H}$, $\sup _{\Gamma(P)} \mathcal{D}=\bigcup \mathcal{D}$.
   \begin{proof}
   Note that each member of $\mathcal{H}$ is a Scott-closed subset of $P$. So to prove the equation, it suffices to show that $\bigcup \mathcal{D} \in \Gamma (P)$. Obviously, $\bigcup \mathcal{D}$ is a lower set. Now Let $D$ be any directed subset of $P$ contained in $\bigcup \mathcal{D}$ such that $\sup D$ exists in $P$. We want to prove that $\sup D \in C$ for some $C \in \mathcal{D}$. First note that $\mathcal{C}=\{ \da d~|~ d\in D\}$ is a directed subset of $\Gamma (P)$. Moreover, $\mathcal{C} \subseteq \mathcal{D}$ by the definition of $  \mathcal{D}$. Since $\mathcal{D}$ is  a $d$-closed set, so $\sup \mathcal{C} \in \mathcal{D}$. But $\sup \mathcal{C}$ is precisely $\da \sup D$. Hence $\sup  D \in C$ for some $C \in \mathcal{D}$.
   \end{proof}
\end{proposition}

\begin{proposition}
   Let $P$ be a poset. If $A \in \mathcal{K}(\Gamma (P))$, then $A$ is an irreducible closed set.
   \begin{proof}
  Suppose $A$ is pre-$C$-compact and $ B,C \in \Gamma(P) $, $A\subseteq B\cup C$. Let $\mathcal{D}=\da \{B,C\}$. Then $\mathcal{D} \in \mathcal{H} $ and $\sup \mathcal{D}=B\cup C$. Hence, $A \in \mathcal{D}$, so $A\subseteq B$ or $A\subseteq C$.
   \end{proof}
\end{proposition}
\begin{lemma}\label{thm}
Let $P$ be a poset and $D$ be a directed subset of $ P$. Then $cl(D)\in \mathcal{K}(\Gamma (P))$.
  \begin{proof}
  Let $ \mathcal{D} \in \mathcal{H}$. Suppose $cl(D)\subseteq \sup  \mathcal{D}$. Then by Proposition \ref{aaa}, $D \subseteq \sup  \mathcal{D}$ and $\{\da d~ |~d\in D \}$ is a directed subset of $\mathcal{D}$. Hence, $cl(D)= \sup  \{ \da d ~| ~d\in D \} \in \mathcal{D}$.
  \end{proof}
\end{lemma}
\begin{corollary}\label{thn}
 Let $P$ be a poset. Then for each $x \in P$, $\da x\prec \da x$ holds.
\end{corollary}
\begin{theorem}\label{bbb}
  For any poset $P$, $\mathcal{K}(\Gamma(P))$ is a dcpo with the inclusion order.
 \begin{proof}
 Let $ (A_i)_{i \in I}$ be a directed subset in $\mathcal{K}(\Gamma(P))$. It suffices to show that $\sup _{i \in I} A_i \prec \sup _{i \in I} A_i$. Suppose $ \mathcal{D} \in \mathcal{H}$ with $\sup _{i \in I} A_i\subseteq \sup  \mathcal{D}$. Then $A_i \subseteq \sup  \mathcal{D}$ for all $ i \in I$. Since $A_i \in \mathcal{K} (\Gamma(P))$, it follows that $A_i \prec A_i$ and so $A_i \in \mathcal{D}$. Because $\mathcal{D}$ is $d$-closed, this implies that $\sup _{i \in I} A_i \in \mathcal{D}$.
 \end{proof}
\end{theorem}

  We write $\hat{P}=\mathcal{K}(\Gamma(P))$ with the relative topology of $\Gamma(P)$.

\begin{proposition}
Let $P$ be a poset. Then $\hat{P}$ is a monotone convergence space.
\begin{proof}
  It is easy to prove by Proposition~\ref{iii}.
\end{proof}

\end{proposition}
\begin{lemma}
  Let $P$ be a poset. If $A\subseteq P$ is tapered and closed, then $A \in \mathcal{K}(\Gamma(P))$.
  \begin{proof}
  Assume $P^\vee$ is the set of all tapered closed subsets and $A \in P^\vee$, then $(P^\vee ,j)$ is the $D$-completion of $P$ by Lemma~\ref{jjj}. Define $f:P\rightarrow \mathcal{K}(\Gamma(P))$ by $f(x)= \da x$ for all $x \in P$. Then, immediately, $f$ is continuous. As $(P^\vee ,j)$ is the $D$-completion of $P$, it follows that there exists a unique continuous map $\hat{f}:P^\vee \rightarrow \mathcal{K}(\Gamma(P))$ s.t. $f=\hat{f} \circ j$. Thus sup$ f(A)=\hat{f}(A)$ by Theorem 3.10 (2) of \cite{Ho03}. We need to prove $A\in \mathcal{K}(\Gamma(P))$. Suppose not, $A \notin \mathcal{K}(\Gamma(P))$. Since $\da a \leq A$ for any $a\in A$ and $\hat{f}$ is continuous, we have $\da a=\hat{f}(\da a)\leq \hat{f}(A)=\sup f(A)$. Hence, $A\varsubsetneqq \sup  f(A)$, i.e. $\sup  f(A) \cap A^c \neq {\O}$. So we have $\sup  f(A) \in \diamondsuit(A^c)$. Because $\hat{f}$ is continuous, this implies $\hat{f}^{-1}(\diamondsuit(A^c))$ is open in $P^\vee$ and $A \in \hat{f}^{-1}(\diamondsuit(A^c))$. Hence, there exists $U \in \sigma (P)$ such that $A \in \diamondsuit U \subseteq \hat{f}^{-1}(\diamondsuit(A^c))$. This means that $A\cap U \neq \emptyset $, then there exists $a \in A \cap U$ such that $\da a \in \diamondsuit U \subseteq \hat{f}^{-1}(\diamondsuit(A^c))$, that is $\hat{f}(\da a )=\da a\in \diamondsuit(A^c)$ which contradicts $a \in A$. Thus we have shown the Lemma.
   \end{proof}
   \end{lemma}
  \begin{theorem}\label{thp}
   Let $P$ be a poset. Then $\mathcal{K}(\Gamma(P))=cl_d(\eta(P))$ where $\eta:P\rightarrow \Gamma (P)$ is defined by $\eta(x)=\da x$ for all $x \in P$,  i.e. $(\hat{P},\eta )$ is the $D$-completion of $P$.
 \begin{proof}
  For any $x \in P$, $\eta(x)=\da x \in \mathcal{K}(\Gamma(P))$ by Corollary \ref{thn}.  Then $\eta(P) \subseteq \mathcal{K}(\Gamma(P))$. We have shown that $\mathcal{K}(\Gamma(P))$ is a subdcpo of $\Gamma(P)$, so $cl_d(\eta(P))\subseteq \mathcal{K}(\Gamma(P))$.

  Conversely, note that $cl_d(\eta(P))\in \mathcal{H}$ and $\sup cl_d(\eta(P))=P$. If $A \in \mathcal{K}(\Gamma(P))$, then $A\subseteq \sup  cl_d(\eta(P))$. This implies $A \in cl_{d}(\eta(P))$ since $A \prec A$.
  \end{proof}
  \end{theorem}
From the proof of Theorem~\ref{thp}, we have the following:
\begin{corollary}\label{coro1}
Let $P$ be a poset. Then the pre-C-compact elements are exactly the tapered closed sets.
\end{corollary}
\begin{definition}
Let $(X, \tau)$ be a topological space. $X$ is called \emph{determined by its directed subsets} if $\tau =\{ U\subseteq X~\mid cl(D)\cap U\neq \emptyset ~\mathrm{implies}~D\cap U \neq \emptyset \}$ for all directed subsets with respect to the order of specialization.
\end{definition}
The following two theorems have the similar proof to that of posets.
\begin{theorem}
Let $X$ be a topological space. Suppose $X$ is determined by its directed subsets, then we can get $(\mathcal{K}(\Gamma(X)),\eta)$ is the $D$-completion of $X$, where $\eta:X\rightarrow \Gamma (X)$ is defined by $\eta(x)=\da x$ for all $x \in X$.
\end{theorem}

\begin{theorem}\label{d}
Suppose $X$ is a topological space. Define $\mathcal{H} =\{\mathcal{D} \subseteq \Gamma(X) \mid \forall A\in \mathcal{D}, \forall a \in A ,\da a \in \mathcal{D}, cl_{d}(\mathcal{D})=\mathcal{D},cl({\bigcup \mathcal{D}})=\bigcup \mathcal{D} \}$. Then $(\mathcal{K}(\Gamma(X)),\eta)$ is the $D$-completion of $X$.
\end{theorem}
\begin{corollary}
  A topological space $X$ is a monotone convergence space if and only if for any $ A \in \mathcal{K}(\Gamma(X))$, there exists a unique point $x$ in $X$, such that $A= \da x$.
\end{corollary}

\section{Reference}
\bibliographystyle{plain}

\end{document}